\documentclass[11pt]{amsart}
\usepackage{amssymb,latexsym}
\usepackage{enumerate}
\usepackage[T1]{fontenc}
\usepackage{mathrsfs}
\usepackage{amsmath}
\usepackage{yhmath}
\makeatletter
\@namedef{subjclassname@2010}{%
  \textup{2010} Mathematics Subject Classification}
\makeatother
\makeindex

\usepackage{chngcntr}

\usepackage{mathtools}
\counterwithin{equation}{section}

\newtheorem{theorem}{Theorem}[subsection]

\theoremstyle{definition}

\newtheorem{sit}[theorem]{Situation}
\newtheorem{rem}[theorem]{Remark}

\newtheorem{prob}[theorem]{Problem}

\numberwithin{theorem}{section}

\newcommand{\mb}{\mathbb}
\newcommand{\mc}{\mathcal}

\newcommand{\s}{\subset}

\begin{document}
\title{Families of strictly pseudoconvex domains and peak functions}
\author{Arkadiusz Lewandowski}
\address{Institute of Mathematics\\ Faculty of Mathematics and Computer Science\\ Jagiellonian University\\ {\L}ojasiewicza 6,
30-348 Kraków, Poland}
\email{Arkadiusz.Lewandowski@im.uj.edu.pl}
\begin{abstract}
We prove that given a family $(G_t)$ of strictly pseudoconvex domains  varying in $\mc{C}^2$ topology on domains, there exists a continuously varying family of peak functions $h_{t,\zeta}$ for all $G_t$ at every $\zeta\in\partial G_t.$  
\end{abstract}



\maketitle
\section{Introduction} Let $D\s\mb{C}^n$ be a bounded domain and let $\zeta$ be a boundary point of $D$. It is called a \emph{peak point} with respect to $\mc{O}(\overline{D})$, the family of functions which are holomorphic in a neighborhood of $\overline{D},$ if there exist a function $f\in\mc{O}(\overline{D})$ such that $f(\zeta)=1$ and $f(\overline{D}\setminus\{\zeta\})\s\mb{D}:=\{z\in\mb{C}:|z|<1\}.$ Such a function is a \emph{peak function for $D$ at $\zeta$}. The concept of peak functions appears to be a powerful tool in complex analysis with many applications. It has been used to show the existence of (complete) proper holomorphic embeddings of strictly pseudoconvex domains into the unit ball $\mb{B}^N$ with large $N$ (see [\ref{For}],[\ref{DD}]), to estimate the boundary behavior of Carath\'{e}odory and Kobayashi metrics ([\ref{Bed}],[\ref{Gra}]), or to construct the solution operators for $\overline{\partial}$ problem with $L^{\infty}$ or H\"{o}lder estimates ([\ref{Forn}],[\ref{Ran2}]), just to name a few of those applications.\\
\indent It is well known that every boundary point of strictly pseudoconvex domain is a peak point. Even more is true, in [\ref{Gra}] it is showed that, given a strictly pseudoconvex domain $G$, there exists an open neighborhood $\widehat{G}$ of $G$, and a continuous function $h:\widehat{G}\times\partial G\rightarrow\mb{C}$ such that for $\zeta\in\partial G$, the function $h(\cdot;\zeta)$ is a peak function for $G$ at $\zeta$.\\
\indent In a recent paper [\ref{DGZ}] the following question has been posed: 
\begin{prob} Let $\rho:\mb{D}\times\mb{C}^n\rightarrow\mb{R}$ be a plurisubharmonic function of class $\mc{C}^{2+k}, k\in\mb{N}\cup\{0\},$ such that for any $z\in\mb{D}$ the truncated function $\rho|_{\{z\}\times\mb{C}^n}$ is strictly plurisubharmonic. Define $G_z:=\{w\in\mb{C}^n:\rho(z,w)<0\},z\in\mb{D}.$ This can be understood as a family of strictly pseudoconvex domains over $\mb{D}.$ 
Does there exist a $\mc{C}^k$-continuously varying family $(h_{z,\zeta})_{z\in\mb{D},\zeta\in\partial G_z}$ of peak functions for $G_z$ at $\zeta$? 
\label{Problem}
\end{prob}
We answer this question affirmatively in the case $k=0$ and under additional assumption that, rougly speaking, the function $\rho$ keeps its regularity up to the set $\Omega\times\mb{C}^n$, where $\Omega$ is some open neighborhood of $\overline{\mb{D}}.$ Namely, let us consider the following:
\begin{sit}
Let $(G_t)_{t\in T}$ be a family of bounded strictly pseudoconvex domains, where $T\s\mb{C}$ is a compact set. Suppose we have a domain $U\s\s\mb{C}^n$ such that
\begin{enumerate}
\item[{(1)}] $\displaystyle{\bigcup_{t\in T}\partial G_t\s\s U},$
\item[{(2)}] for each $t\in T$ there exists a defining function $r_t$ for $G_t$ satisfying with neighborhood $\partial G_t\s U$ all the conditions {(A)-(F)} below (see Section 2),
\item[{(3)}] for any $\varepsilon>0$ there exists a $\delta>0$ such that for any $s,t\in T$ with $|s-t|\leq \delta$ there is $\|r_t-r_s\|_{\mc{C}^2(U)}<\varepsilon$.\label{C2}
\end{enumerate}\label{Situation}
\end{sit}
Observe that the above setting is completely in the spirit of the formulation of Problem \ref{Problem}:
\begin{enumerate}[(i)]
\item The assumption that all the functions $r_t$ satisfy (A)-(F) with common neighborhood $\partial G_t\s U$ stays in relation with the fact that in Problem \ref{Problem} all the defining functions for domains $D_z$ have the same domain of definition ($\mb{C}^n$).
\item The assumption (3) comes from the fact that the function $\rho$ in Problem \ref{Problem} is of class at least $\mc{C}^2.$
\item The compactness of the set of parameters ($T$) reflects the above mentioned assumption that $\rho$ continues to be of class $\mc{C}^2$ up to $\Omega\times\mb{C}^n$, with $\Omega$ being some neighborhood of $\overline{\mb{D}}.$  
\end{enumerate}
We shall prove the following:
\begin{theorem}
Let $(G_t)_{t\in T}$ be a family of strictly pseudoconvex domains as in Situation \ref{Situation}.
Then there exists an $\varepsilon>0$ such that for any $\eta_1<\varepsilon$ there exist an $\eta_2>0$ and positive constants $d_1,d_2$ such that for any $t\in T$ there exist a domain $\widehat{G_t}$ containing $\overline{G_t}$, and functions $h_t(\cdot;\zeta)\in\mc{O}(\widehat{G_t}),\zeta\in\partial G_t$ fulfilling the following conditions:
\begin{enumerate}
{\item[\emph{(a)}] $h_t(\zeta;\zeta)=1, |h_t(\cdot;\zeta)|<1$ on $\overline{G_t}\setminus\{\zeta\}$ (in particular, $h_t(\cdot;\zeta)$ is a peak function for $G_t$ at $\zeta$), \label{a}}
{\item[\emph{(b)}] $|1-h_t(z;\zeta)|\leq d_1\|z-\zeta\|, z\in\widehat{G_t}\cap\mb{B}(\zeta,\eta_2),$\label{b}}
{\item[\emph{(c)}] $|h_t(z;\zeta)|\leq d_2<1, z\in\overline{G_t},\|z-\zeta\|\geq\eta_1.$\label{c}}
\end{enumerate}
Moreover, the constants $\varepsilon, \eta_2, d_1, d_2$, domains $\widehat{G_t},$ and functions $h_t(\cdot;\zeta)$ may be chosen in such a way that for any $\alpha>0$ and any fixed triple $(t_0,\zeta_0,z_0)$, where $t_0\in T, \zeta_0\in\partial G_{t_0},$ and $z_0\in\widehat{G_{t_0}}$, there exists a $\delta>0$ such that whenever the triple $(s, \xi, w)$ satisfies $s\in T, \xi\in\partial G_s, w\in\widehat{G_s},$ and $\max\{|s-t_0|,\|\xi-\zeta_0\|,\|w-z_0\|\}<\delta$, then
$|h_{t_0}(z_0;\zeta_0)-h_s(w;\xi)|<\alpha$.
\label{Main}
\end{theorem}
The latter property will be referred to as \emph{continuity}.
\begin{rem}
It is known that for each $t\in T$ there exists an $\varepsilon=\varepsilon(t)>0$ such that for any $\eta_1<\varepsilon$ there exist a positive $\eta_2=\eta_2(t)<\eta_1$, constants $d_1=d_1(t), d_2=d_2(t)\in\mb{R},$ domain $\widetilde{G_t}$ containing $\overline{G_t}$, and functions $h_t(\cdot;\zeta)\in\mc{O}(\widehat{G_t}),\zeta\in\partial G_t$ satisfying (a)-(c). This is a subject of Theorem  19.1.2 from [\ref{JP}]. The strength of our result dwells in the fact that all the constants $\varepsilon, \eta_2, d_1, d_2$ are chosen independently of $t$ and in the continuity property.
\end{rem}
In Section 2 we recall some preliminaries concerning the strictly pseudoconvex domains. The proof of Theorem \ref{Main} is presented in Section 3.
\section{Strictly pseudoconvex domains} 
Let $D\s\s\mb{C}^n$ be a domain. It is called a \emph{strictly pseudoconvex} if there exist a neighborhood $U$ of $\partial D$ and a \emph{defining function} $r:U\rightarrow\mb{R}$ of class $\mc{C}^2$ and such that
\begin{enumerate}
{\item[(A)] $D\cap U=\{z\in U:r(z)<0\}$,\label{Condition1}}
{\item[(B)] $(\mb{C}^n\setminus\overline{D})\cap U=\{z\in U:r(z)>0\},\label{Condition2}$}
{\item[(C)] $\nabla r(z)\neq 0$ for $z\in \partial D,$ where $\nabla r(z):=\left(\frac{\partial r}{\partial\overline{z}_1}(z),\cdots,\frac{\partial r}{\partial\overline{z}_n}(z)\right)$,\label{Condition3}}
\end{enumerate}
together with $$\mc{L}_r(z;X)>0\text{\ for\ } z\in\partial D\text{\ and\ nonzero\ }X\in T_z^{\mb{C}}(\partial D),$$  
where $\mc{L}_r$ denotes the Levi form of $r$ and $T_z^{\mb{C}}(\partial D)$ is the complex tangent space to $\partial D$ at $z$.\\ 
\indent It is known that $U$ and $r$ can be chosen to satisfy (A)-(C) and, additionally:
\begin{enumerate}
{\item[(D)] $\mc{L}_r(z;X)>0$ for $z\in U$ and all nonzero $X\in\mb{C}^n,$\label{Condition4}}
{\item[(E)] $\| \nabla r(z)\|=1,z\in\partial D,$ \label{Condition5}}
{\item[(F)] for every $z\in D\cap U$ there is a unique $\pi(z)\in\partial D$ with $$\text{dist}(z,\partial D)=\|z-\pi(z)\|,\label{Condition6}$$}
\end{enumerate}
cf. [\ref{Kra1}],[\ref{Kra2}]. Note that for a function $r$ as above and a point $\zeta\in\partial G$, Taylor expansion of $r$ at $\zeta$ has the following form:
\begin{equation}
r(z)=r(\zeta)-2\text{Re}P(z;\zeta)+\mc{L}_{r}(\zeta;z-\zeta)+o(\|z-\zeta\|^2),\label{Taylor}
\end{equation}
where 
$$\displaystyle{
P(z;\zeta):=-\sum_{j=1}^n\frac{\partial r}{\partial z_j}(\zeta)(z_j-\zeta_j)-\frac{1}{2}\sum_{i,j=1}^n\frac{\partial^2r}{\partial z_i\partial z_j}(\zeta)(z_i-\zeta_i)(z_j-\zeta_j)
}$$
is the \emph{Levi polynomial of $r$ at $\zeta$}.
\section{Proof of Theorem \ref{Main}}
We divide the proof into two parts. First we give the construction of $\widehat{G_t}$ and $h_t(\cdot;\zeta),t\in T$, and define the constants $\varepsilon, \eta_2, d_1$, and $d_2$, all independent of $t$. This is refinement of the construction from the proof of Theorem 19.1.2 from [\ref{JP}]. Note that in order to get the independence of all the constants from $t$, we must be more careful here. In the second part we prove the continuity property.
\begin{proof}[Construction of $\widehat{G_t}$ and $h_t(\cdot;\zeta)$ and the choice of $\varepsilon, \eta_2, d_1$, and $d_2$] For $t\in T$ and $\zeta\in\partial G_t$ let $P_t(z;\zeta)$ be the Levi polynomial of $r_t$ at $\zeta$.\\
\indent Fix an $\varepsilon_1>0$ such that $\displaystyle{U':=\bigcup_{t\in T,\zeta\in\partial G_t}\mb{B}(\zeta,\varepsilon_1)}\s\s U.$\\
There exists a constant $C_1=C_1(t)<1$ such that
$$
\mc{L}_{r_t}(z;X)\geq C_1\|X\|^2, \quad z\in U',X\in\mb{C}^n.
$$
Indeed, $\mc{L}_{r_t}$ is continuous and positive on $U\times(\mb{C}^n\setminus\{0\})$, so it attains its minimum $C_1(t)>0$ on $\overline{U'}\times\mb{S}^{n-1}.$ Since for any nonzero $X\in\mb{C}^n$ we have $\frac{X}{\|X\|}\in\mb{S}^{n-1},$ we get the required inequality. Moreover, from the assumption (3) it follows that for $s$ from some neighborhood of $t$ we have
$$
\mc{L}_{r_s}(z;X)\geq\frac{C_1(t)}{2}\|X\|^2, \quad z\in U',X\in\mb{C}^n.
$$
The compactness argument then gives that $C_1$ may be chosen independently of $t$.\\
\indent Taylor formula (\ref{Taylor}) yields that with some $0<C_2<C_1$ there is
\begin{equation}
r_t(z)\geq -2\text{Re}P_t(z;\zeta)+C_2\|z-\zeta\|^2 \label{EC2}
\end{equation}
for $\|z-\zeta\|<\varepsilon_2(t)<\varepsilon_1,\zeta\in\partial G_t$, where $\varepsilon_2(t)$ is independent of $\zeta\in\partial G_t$ (and even of $\zeta\in W\s\s U,$ some neighborhood of $\partial G_t$ - see [\ref{Ran}], Proposition II.2.16). Moreover, from the proof of Theorem V.3.6 from $[\ref{Ran}]$ it follows that for $s$ close enough to $t$ we have
$$
r_s(z)\geq r_s(\zeta)-2\text{Re}P_s(z;\zeta)+\frac{C_2}{2}\|z-\zeta\|^2,\quad\zeta\in W,\|z-\zeta\|<\varepsilon_2(t).
$$
Therefore, for $s$ near to $t$, and for $\xi\in\partial G_s$, the following estimate holds true:
$$
r_s(z)\geq - 2\text{Re}P_s(z;\xi)+\frac{C_2}{2}\|z-\xi\|^2,\quad \|z-\xi\|<\varepsilon_2(t).
$$
The compactness argument then implies that $C_2$ and $\varepsilon_2$ in (\ref{EC2}) may be chosen independently of $t$.\\
\indent Let $0<\eta_1<\varepsilon_2$ and $\widehat{\chi}\in\mc{C}^{\infty}(\mb{R},[0,1])$ be such that $\widehat{\chi}(t)=1$ for $t\leq\frac{\eta_1}{2}$ and $\widehat{\chi}(t)=0$ for $t\geq\eta_1.$ Put 
$\chi(z;\zeta):=\widehat{\chi}(\|z-\zeta\|).$ This is a smooth function on $\mb{C}^n\times\mb{C}^n$, taking its values in $[0,1]$.\\
\indent Define
$$
\varphi_t(z;\zeta):=\chi(z;\zeta)P_t(z;\zeta)+(1-\chi(z;\zeta))\|z-\zeta\|^2,\quad z\in\mb{C}^n.
$$
Observe that if $\|z-\zeta\|\leq\frac{\eta_1}{2},$ then $\varphi_t(z;\zeta)=P_t(z;\zeta).$ In particular $\varphi_t(\cdot;\zeta)\in\mc{O}(\mb{B}(\zeta,\frac{\eta_1}{2}))$. Furthermore, for $z$ satisfying $\|z-\zeta\|\geq\frac{\eta_1}{2}$ and $r_t(z)< C_2\frac{\eta_1^2}{8}$ the following estimate holds true:
\begin{equation}
2\text{Re}\varphi_t(z;\zeta)\geq C_2\frac{\eta_1^2}{8}>0. \label{2Re}
\end{equation}
Take $0<\eta_t<C_2\frac{\eta_1^2}{8}$ such that the connectend component $\widetilde{G_t}$ containing $\overline{G_t}$ of the open set
$$
G_t\cup\{z\in U':r_t(z)<\eta_t\}
$$ 
is a strictly pseudoconvex domain, relatively compact in $G_t\cup U'.$ Because of the assumption (3), there exists a positive number $\beta$ such that for $s$ close to $t$ the connected component $\widetilde{G_s}$ containing $\overline{G_s}$ of the set
$$
G_s\cup\{z\in U':r_s(z)<\eta_t-\beta\}
$$
is a strictly pseudoconvex domain, relatively compact in $G_s\cup U'.$ Making again use of the compactness of $T$, we conclude that in fact $\eta=\eta_t$ may be taken independently of $t.$ Note that for the family $(\widetilde{G_t})_{t\in T}$ the assumption (3) remains true.\\
\indent The function $\varphi_t(\cdot;\zeta)\in\mc{C}^{\infty}(\mb{C}^n)$ does not vanish on $\widetilde{G}\setminus\mb{B}(\zeta,\frac{\eta_1}{2})$ and is in $\mc{O}(\mb{B}(\zeta,\frac{\eta_1}{2}))$. Therefore $\bar{\partial}\frac{1}{\varphi_t(\cdot;\zeta)}$ defines a $\bar{\partial}$-closed $\mc{C}^{\infty}$ form $$\displaystyle{\alpha_t(\cdot;\zeta)=\sum_{j=1}^n\alpha_{t,j}(\cdot;\zeta)dz_j}$$ on $\widetilde{G_t},$ where
$$\displaystyle{
\alpha_{t,j}=\begin{cases}0,&z\in\widetilde{G_t}\cap\mb{B}(\zeta;\frac{\eta_1}{2}),\\
-\frac{\partial\varphi_t}{\partial\bar{z_j}}(z;\zeta)\cdot\frac{1}{\varphi_t^2(z;\zeta)},&z\in\widetilde{G_t}\setminus\mb{B}(\zeta;\frac{\eta_1}{2}).
\end{cases}
}$$
Thanks to (\ref{2Re}) we have $\|\alpha_{t,j}(\cdot;\zeta)\|_{\widetilde{G_t}}\leq C_3$, where, utilizing the compactness of $T$ together with the assumption (3), we deliver that $C_3$ is independent of $t$ and $\zeta\in\partial G_t$.
[\ref{Ran}, Theorem V.2.7] gives then the functions $v_t(\cdot;\zeta)\in\mc{C}^{\infty}(\widetilde{G_t})$ with $\bar{\partial}v_t(\cdot;\zeta)=\alpha_t(\cdot;\zeta)$
and
$$
\|v_t(\cdot;\zeta)\|_{\widetilde{G_t}}\leq C_4,
$$
where $C_4$ does not depend on $\zeta\in\partial G_t.$ Moreover, by [\ref{Ran}, Theorem V.3.6] and the compactness of $T$, $C_4$ may be chosen to be independent of $t$.\\
Define
$$
f_t(\cdot;\zeta):=\frac{1}{\varphi_t}(\cdot;\zeta)+C_4-v_t(\cdot;\zeta),\quad z\in\widetilde{G_t}\setminus Z_t(\zeta),
$$
where
$$
Z_t(\zeta):=\{z\in\widetilde{G_t}:\varphi_t(z;\zeta)=0\}.
$$
Then $f_t(\cdot;\zeta)\in\mc{O}(\widetilde{G_t}\setminus Z_t(\zeta))$ as well as
$$
\text{Re}f_t(\cdot;\zeta)>0
$$
on the set $(\widetilde{G_t}\setminus\mb{B}(\zeta,\frac{\eta_1}{2}))\cup(\overline{G_t}\setminus\{\zeta\})$, in virtue of (\ref{EC2}) and (\ref{2Re}). Since for any $\zeta\neq z_0\in\partial{G_t}\cap\overline{\mb{B}}(\zeta,\frac{\eta_1}{2})$ there exists a neighborhood $U_{z_0}$ of $z_0$ such that $\text{Re}f_t(\cdot;\zeta)>0$ on $U_{z_0}$, we conclude that there exists a neighborhhod $U_{t,\zeta}$ of $\overline{G_t}\setminus\{\zeta\}$ such that the function
$$
h_t(\cdot;\zeta):=\text{exp}(-g_t(\cdot;\zeta)),
$$
where $g_t(\cdot;\zeta):=\frac{1}{f_t(\cdot;\zeta)}$, is holomorphic on $H_{t,\zeta}:=(\widetilde{G_t}\setminus\mb{B}(\zeta,\frac{\eta_1}{2}))\cup U_{t,\zeta}.$ Note that $h_t$ takes its values in $\mb{D}.$\\
\indent There exists a $C_5>0,$ independent on $t$, such that
$$
|P_t(z;\zeta)|\leq C_5\|z-\zeta\|,\quad \zeta, z\in U'.
$$
Therefore, since for $0<\eta_2<\min\big\{\frac{\eta_1}{2},\frac{1}{4C_4C_5}\big\}$, which now is independent of $t$, and for $z\in ((H_{t,\zeta}\cup\mb{B}(\zeta,\eta_2))\cap\mb{B}(\zeta,\eta_2))\setminus Z_t(\zeta)$
the following equality holds true:
$$
g_t(z;\zeta)=\frac{P_t(z;\zeta)}{1-P_t(z;\zeta)(v_t(z;\zeta)-C_4)},
$$
we conclude that $g_t(\cdot;\zeta)$ is bounded near $Z_t(\zeta)$, which yields it extends to be holomorphic on $\widetilde{H_{t,\zeta}}:=H_{t,\zeta}\cup\mb{B}(\zeta;\eta_2).$\\
\indent Now $\widetilde{H_{t,\zeta}}$ depends on $\zeta$, but using the inclusion $\overline{G_t}\s\widetilde{H_{t,\zeta}}$, we may find some $\widehat{G_t},$ strictly pseudoconvex domain which is independent on $\zeta\in\partial G_t$, such that $\overline{G_t}\s\widehat{G_t}\s\widetilde{H_{t,\zeta}}$ for each $\zeta\in\partial G_t$, and with the property that $h_t(\cdot;\zeta)\in\mc{O}(\widehat{G_t}), \zeta\in\partial G_t$ (use the joint continuity of $\varphi_t$ with respect to $z$ and $\zeta$ to shrink $\widetilde{H_{t,\zeta}}$ little bit to get some domain with desired properties, independent on $\xi$ close to $\zeta$, and finally apply the compactness of $\partial G_t$).\\
Let $C_6$, independent on $t$ and $\zeta\in\partial G_t$, such that for $z\in\widehat{G_t}$ with $\|z-\zeta\|<\eta_2$ we have
$$
g_t(z;\zeta)\leq\frac{C_5\|z-\zeta\|}{1-2C_4C_5\|z-\zeta\|}\leq C_6\|z-\zeta\|.
$$
This implies
$$
|1-h_t(z;\zeta)|\leq C_7|g_t(z;\zeta)|\leq C_6 C_7\|z-\zeta\|=:d_1\|z-\zeta\|
$$
for $z\in\widehat{G_t},\|z-\zeta\|<\eta_2,\zeta\in\partial G_t$, if only $C_7$ is chosen so that
$$
|e^{\lambda}-1|\leq C_7|\lambda|,\quad |\lambda|\leq C_6\eta_2.
$$
In particular, $d_1$ does not depend on $t$ and we have $h_t(\zeta;\zeta)=1.$\\
\indent Furthermore, for $z\in\overline{G_t},\|z-\zeta\|\geq\eta_1$ there is
\begin{multline*}
\text{Re}g_t(z;\zeta)=\|z-\zeta\|^2\frac{1+\|z-\zeta\|^2(C_4-\text{Re}v_t(z;\zeta))}{|1-\|z-\zeta\|^2(v_t(z;\zeta)-C_4)|^2}\\
\geq\frac{\eta_1^2}{(1+2(\text{diam}U)^2C_4)^2}=:C_8,
\end{multline*}
which gives 
\begin{equation}
|h_t(z;\zeta)|\leq e^{-C_8}=:d_2<1.\label{D2}
\end{equation}
Observe that $d_2$ is independent on $t.$
\end{proof}
\begin{proof}[Proof of continuity] Fix $\alpha>0, t_0\in T,\zeta_0\in\partial G_{t_0},$ and $z_0\in\widehat{G_{t_0}}$. Let $K_0$ be a compact subset of $\widetilde{G_{t_0}},$ containing in its interior the set $\overline{G_{t_0}}\cup\{z_0\}.$ In the sequel we shall use the following convention: whenever we say that the triple $(s,\xi,w)$ is near to $(t_0,\zeta_0,z_0),$ it will carry the additional information that $\xi\in\partial G_s,w\in\widehat{G_s}$, unless explicitly stated otherwise.\\
Observe that for $(s,\xi)$ close to $(t_0,\zeta_0)$ (even without requiring that $\xi\in\partial G_s$), and any $z,w\in U'$ we have
$$
|P_{t_0}(z;\zeta_0)-P_s(w;\xi)|<M_1\alpha
$$
with some positive $M_1.$ In particular, the same estimate is true for $z=z_0$ and $w$ close to $z_0.$\\ 
Further, using the fact that all the functions $\varphi_t$ are continuous as functions of both variables, we conclude that for $(s,\xi)$ close to $(t_0,\zeta_0)$ we have 
$$
\|\varphi_{t_0}(\cdot;\zeta_0)-\varphi_s(\cdot;\xi)\|_{U'}<M_2\alpha
$$
with some positive $M_2.$\\
\indent For $(s,\xi)$ near $(t_0,\zeta_0)$ we have
$$
\big\|\frac{\partial\varphi_{t_0}}{\partial\bar{z_j}}(\cdot;\zeta_0)-\frac{\partial\varphi_{s}}{\partial\bar{z_j}}(\cdot;\xi)\big\|_{U'}<M_3\alpha
$$
with some positive $M_3.$ Furthermore, for $z\in\widetilde{G_s}\cap\widetilde{G_t}$ the following estimates hold true:\\

(I) If $z\notin\mb{B}(\zeta_0,\frac{\eta_1}{2})\cup\mb{B}(\xi,\frac{\eta_1}{2})$, then
$$
|\alpha_{t_0,j}(z;\zeta_0)-\alpha_{s,j}(z;\xi)|< L\alpha,
$$
where positive constant $L$ does not depend on $z$ as above. Indeed,
\begin{multline*}
\left|\frac{\frac{\partial\varphi_{t_0}}{\partial\bar{z_j}}(z;\zeta_0)}{\varphi^2_{t_0}(z;\zeta_0)}-\frac{\frac{\partial\varphi_{s}}{\partial\bar{z_j}}(z;\xi)}{\varphi^2_{s}(z;\xi)}\right|=
\left|\frac{\varphi_s^2(z;\xi)\frac{\partial\varphi_{t_0}}{\partial\bar{z_j}}(z;\zeta_0)-\varphi_{t_0}^2(z;\zeta_0)\frac{\partial\varphi_{s}}{\partial\bar{z_j}}(z;\xi)}{\varphi^2_{t_0}(z;\zeta_0)\varphi_s^2(z;\xi)}\right|\\
\leq\frac{64}{C_2^2\eta_1^4}\big|\varphi_s^2(z;\xi)\frac{\partial\varphi_{t_0}}{\partial\bar{z_j}}(z;\zeta_0)-\varphi_{t_0}^2(z;\zeta_0)\frac{\partial\varphi_{s}}{\partial\bar{z_j}}(z;\xi)\big|\\
\leq\frac{64}{C_2^2\eta_1^4}\|\varphi^2_s\|_{U'}\big\|\frac{\partial\varphi_{t_0}}{\partial\bar{z_j}}(\cdot;\zeta_0)-\frac{\partial\varphi_{s}}{\partial\bar{z_j}}(\cdot;\xi)\big\|_{U'}+\big\|\frac{\partial\varphi_s}{\partial\bar{z_j}}(\cdot;\xi)\big\|_{U'}\|\varphi^2_s-\varphi^2_{t_0}\|\\
\leq \frac{64}{C_2^2\eta_1^4}L_1M_3\alpha+L_2M_2\alpha=:L\alpha,
\end{multline*}
where the first inequality is the consequence of (\ref{2Re}).\\

(II) If $z\in\mb{B}(\zeta_0,\frac{\eta_1}{2})\cup\mb{B}(\xi,\frac{\eta_1}{2})$:\\ Observe that letting $\xi$ close to $\zeta_0$, we may make the balls arbitrarily close each other. Using then the assumption (3), the fact that $\eta$ were chosen to be strictly smaller than $C_2\frac{\eta_1^2}{8}$, and the strictness of uniform estimate (\ref{2Re}), we see that for $(s,\xi)$ close enough to $(t_0,\zeta_0)$
the estimate similar to the previous one holds true for $\displaystyle{z\in S:=\bigcup_{w:\|w-\zeta_0\|=\frac{\eta_1}{2}}\mb{B}(w,\gamma)}$ with some sufficiently small $\gamma>0$, (and is independent on such $z$). Additionally, $(s,\xi)$ may be chosen so that $S':=(\mb{B}(\zeta_0,\frac{\eta_1}{2})\cup\mb{B}(\xi,\frac{\eta_1}{2}))\setminus S\s \mb{B}(\zeta_0,\frac{\eta_1}{2})\cap\mb{B}(\xi,\frac{\eta_1}{2})$.\\ 

Noting that for $z\in S'$ and $(s,\xi)$ as above $\alpha_{t_0,j}(z;\zeta_0)=\alpha_{s,j}(z;\xi)=0,$ we conclude that
$$
\|\alpha_{t_0}(\cdot;\zeta_0)-\alpha_s(\cdot;\xi)\|_{\widetilde{G_{t_0}}\cap\widetilde{G_s}}\leq M_4\alpha
$$
with some positive $M_4.$\\

Ofcourse $\overline{G_{t_0}}\s\widetilde{G_{t_0}}$. This yields that for $s$ close to $t_0$ we have $\overline{G_{t_0}}\s\widetilde{G_s}$ as well as $\overline{G_s}\s\widetilde{G_{t_0}}$ (the assumption (3) remains true for the family $(\widetilde{G_t})_{t\in T}$).  For $s$ close to $t_0$ we may now pick some $G_{t_0,s}$, a strictly pseudoconvex domain with smooth boundary and such that
$$
\overline{G_s}\cup\overline{G_{t_0}}\s K_{0}\s\s G_{t_0,s}\s\s\widetilde{G_s}\cap\widetilde{G_{t_0}}.
$$
Again thanks to the property (3), $G_{t_0,s}$ may be chosen independently of $s$ if $s$ is close enough to $t_0$. For such $s$, denote it by $G^{t_0}.$ Then, using Lemma 2 from [\ref{Gra}], we find some positive constant $\Gamma$ such that
$$
\|v_{t_0}(\cdot;\zeta_0)-v_s(\cdot;\xi)\|_{K_{0}}\leq\Gamma\|\alpha_{t_0}(\cdot;\zeta_0)-\alpha_s(\cdot;\xi)\|_{G^{t_0}}\leq\Gamma M_4\alpha=:M_5\alpha.
$$ 
Consequently, for $(s,\xi,w)$ close to $(t_0,\zeta_0,z_0)$ there is
$$
|v_{t_0}(z_0;\zeta_0)-v_s(w;\xi)|\leq|v_{t_0}(z_0;\zeta_0)-v_{t_0}(w;\zeta_0)|+|v_{t_0}(w;\zeta_0)-v_s(w;\xi)|\leq M_6\alpha
$$
for some positive $M_6$ (use the smoothness of $v_{t_0}(\cdot;\zeta_0)$).\\

There are two cases to be considered:\\

\emph{Case 1.} $z_0\in H_{t_0,\zeta_0}\cap\text{int}K_{0}$.\\
Then $\varphi_{t_0}(z_0;\zeta_0)\neq 0$ and for $(s,\xi,w)$ near $(t_0,\zeta_0,z_0)$ we have $\varphi_s(w;\xi)\neq 0.$ For such $(s,\xi,w)$ we have
\begin{multline*}
|f_{t_0}(z_0;\zeta_0)-f_s(w;\xi)|\leq\big|\frac{1}{\varphi_{t_0}(z_0;\zeta_0)}-\frac{1}{\varphi_s(w;\xi)}\big|+|v_{t_0}(z_0;\zeta_0)-v_s(w;\xi)|\\
\leq \big|\frac{\varphi_s(w;\xi)-\varphi_{t_0}(z_0;\zeta_0)}{\varphi_{t_0}(z_0;\zeta_0)\varphi_s(w;\xi)}\big|+M_6\alpha.
\end{multline*}
Considering the last but one term, its denominator is bounded below by some positive constant for $(s,\xi,w)$ close to $(t_0,\zeta_0,z_0)$, and the counter is estimated from above by $M_2\alpha.$ Thus for $(s,\xi,w)$ close to $(t_0,\zeta_0,z_0)$
$$
|f_{t_0}(z_0;\zeta_0)-f_s(w;\xi)|\leq M_7\alpha
$$
for some positive $M_7.$\\
\indent In our situation the function $g_{t_0}(\cdot;\zeta_0)$ is holomorphic in a neighborhood of $z_0$ and so is $g_s(\cdot;\xi)$ for $(s,\xi)$ close to $(t_0,\zeta_0).$ We conclude that for $(s,\xi,w)$ close to $(t_0,\zeta_0,z_0)$ there is
$$
|g_{t_0}(z_0;\zeta_0)-g_s(w;\xi)|\leq M_8\alpha
$$
for some positive $M_8$, and 
$$
|h_{t_0}(z_0;\zeta_0)-h_s(w;\xi)|=\big|\text{exp}(-g_{t_0}(z_0;\zeta_0))-\text{exp}(-g_s(w;\xi))\big|\leq M_9\alpha
$$
for some positive $M_9.$\\

\emph{Case 2.} $z_0\in(\widetilde{G_{t_0}}\cap\mb{B}(\zeta_0,\eta_2))\cap\text{int}K_{t_0}.$\\

(I) Suppose $\varphi_{t_0}(z_0;\zeta_0)\neq 0.$\\
It is equivalent to $P_{t_0}(z_0;\zeta_0)\neq 0.$ This yields that $P_s(w;\xi)\neq 0$ for $(s,\xi,w)$ close to $(t_0,\zeta_0,z_0)$. Then
\begin{multline*}
|g_{t_0}(z_0;\zeta_0)-g_s(w;\xi)|\\=\big|\frac{P_{t_0}(z_0;\zeta_0)}{1-P_{t_0}(z_0;\zeta_0)(v_{t_0}(z_0;\zeta_0)-C_4)}-\frac{P_{s}(w;\xi)}{1-P_{s}(w;\xi)(v_{s}(w;\xi)-C_4)}\big|\\
\leq N|P_{t_0}(z_0;\zeta_0)-P_s(w;\xi)|+N|P_{t_0}(z_0;\zeta_0)P_s(w;\xi)\|v_{t_0}(z_0;\zeta_0)-v_s(w;\xi)|\\
\leq NM_1\alpha+NM_6\alpha=:M_{10}\alpha,
\end{multline*}
and similarly as in the previous case
$$
|h_{t_0}(z_0;\zeta_0)-h_s(w;\xi)|\leq M_{11}\alpha
$$
with some positive $N,M_{10}$, and $M_{11}.$\\

(II) Suppose $\varphi_{t_0}(z_0;\zeta_0)=0.$\\
This is equivalent to $P_{t_0}(z_0;\zeta_0)=0.$ Then for some positive $\rho$ we have $\mb{B}(z_0,\rho)\s\s K_0\cap\mb{B}(\zeta_0,\eta_2).$ Similarly, for $(s,\xi)$ close to $(t_0,\zeta_0)$ there is $\mb{B}(z_0,\rho)\s\s K_0\cap\mb{B}(\xi,\eta_2).$ Therefore, because of the choice of $d_2$ in (\ref{D2}), for $(s,\xi,w)$ close to $(t_0,\zeta_0,z_0),w\in\mb{B}(z_0,\rho/2)$ there is
\begin{multline*}
|h_{t_0}(w;\zeta_0)-h_s(w;\xi)|\leq |1-h_{t_0}(w;\zeta_0)|+|1-h_s(w;\xi)|\\ \leq d_1(\|w-\zeta_0\|+\|w-\xi\|)\leq2d_1\eta_2.
\end{multline*}
Consequently, since the functions $h_{t_0}(\cdot;\zeta_0)-h_s(\cdot;\xi)$ are holomorphic in suitable neighborhood of $z_0$ for $(s,\xi)$ close to $(t_0,\zeta_0)$, for some positive $\widetilde{\rho}<\rho/2$, for every $x,y\in\mb{B}(z_0,\widetilde{\rho}/2)$ we have
\begin{align}
|h_{t_0}(x;\zeta_0)-h_s(x;\xi)-h_{t_0}(y;\zeta_0)+h_s(y;\xi)|\leq\alpha.\label{AA}
\end{align}
Moreover, $\widetilde{\rho}$ may be chosen so that for $v,w\in\mb{B}(z_0,\widetilde{\rho}/2)$ there is
\begin{align}
|h_{t_0}(v;\zeta_0)-h_{t_0}(w;\zeta_0)|\leq\alpha,\label{Cont}
\end{align}
by continuity of $h_{t_0}(\cdot;\zeta_0).$\\
Fix some $w_0\in\mb{B}(z_0,\widetilde{\rho}/2)$ such that $P_{t_0}(w_0;\zeta_0)\neq 0.$ Then for $(s,\xi)$ near $(t_0,\zeta_0)$, by virtue of the subcase (I), we have $$|h_{t_0}(w_0;\zeta_0)-h_s(w_0;\xi)|\leq\alpha.$$ Finally, for $w\in\mb{B}(z_0,\widetilde{\rho}/2)$ and $(s,\xi)$ close to $(t_0,\zeta_0)$ we have
\begin{multline*}
|h_{t_0}(z_0;\zeta_0)-h_s(w;\xi)|\leq|h_{t_0}(z_0;\zeta_0)-h_{t_0}(w_0;\zeta_0)|+|h_{t_0}(w_0;\zeta_0)-h_s(w_0;\xi)|\\+|h_s(w_0;\xi)-h_s(w;\xi)|\leq\alpha+\alpha+2\alpha=4\alpha,
\end{multline*}
where the last estimate follows from (\ref{AA}) and (\ref{Cont}), which leads us to the conclusion.
\end{proof}

\end{document}